\documentclass[twocolumn]{autart}

\usepackage{graphicx}  
\usepackage{algorithm,algorithmicx,algpseudocode,amsmath,amsfonts,verbatim,mathtools,enumerate,circuitikz,breakcites,accents,hyperref}
\usepackage[caption=false]{subfig}
\allowdisplaybreaks

\newcommand{\argmin}{\mathop{\rm argmin}}

\newcommand{\mnorm}[1]{{\left\vert\kern-0.25ex\left\vert\kern-0.25ex\left\vert #1 
    \right\vert\kern-0.25ex\right\vert\kern-0.25ex\right\vert}}

\newtheorem{definition}{Definition} 
\newtheorem{theorem}{Theorem}
\newtheorem{lemma}{Lemma}

\newtheorem{remark}{Remark}

\newtheorem{assumption}{Assumption}
\newtheorem{game}{Game}

\algnewcommand\algorithmicforeach{\textbf{for each}}
\algdef{S}[FOR]{ForEach}[1]{\algorithmicforeach\ #1\ \algorithmicdo}

\newcommand{\ie}{{\it i.e.}}

\usepackage{cite}
\usetikzlibrary{matrix,arrows}
\usetikzlibrary{shapes}
\usetikzlibrary{calc}
\usetikzlibrary{trees}
\usetikzlibrary{positioning} 
\usetikzlibrary{decorations.markings}
\usetikzlibrary{patterns}
\usetikzlibrary{shapes.misc, positioning}
\tikzset{
  big arrow/.style={
    decoration={markings,mark=at position 1 with {\arrow[scale=2.0, #1]{>}}},
    postaction={decorate},
    shorten >=0.4pt, line width=0.35mm},
  }
  
\begin{document}
\begin{frontmatter}
\runtitle{run title}

\title{Variable Demand and Multi-commodity Flow in Markovian Network Equilibrium}

\author{Yue Yu},
\author{Dan Calderone},
\author{Sarah H. Q. Li}, 
\author{Lillian J. Ratliff},
\author{Beh\c{c}et A\c{c}\i kme\c{s}e}
\address{Oden Institute for Computational Engineering and Sciences, The University of Texas at Austin, Austin, TX, 78712
}
\address{Department of Aeronautics and Astronautics, University of Washington, Seattle, WA, 98195
}
\address{Department of Electrical and Computer Engineering, University of Washington, Seattle, WA, 98195
}

\thanks{Y. Yu is with the Oden Institute for Computational Engineering and Sciences, The University of Texas at Austin, Austin, Texas, 78712 (e-mail: {\tt yueyu@utexas.edu}). S. H. Q. Li, and B. A\c{c}\i kme\c{s}e are with the William E. Boeing Department of Aeronautics \& Astronautics, University of Washington, Seattle, Washington, 98195 (e-mail: {\tt yueyu@uw.edu, sarahli@uw.edu, behcet@uw.edu}). D. Calderone and L. J. Ratliff are with the Department of Electrical \& Computer Engineering, University of Washington, Seattle, Washington, 98195 (e-mail:{\tt djcal@uw.edu}, {\tt ratliffl@uw.edu}). }

\begin{keyword}                           
Wardrop equilibrium, Markov decision process, network optimization     
\end{keyword}  

\begin{abstract}      
Markovian network equilibrium generalizes the classical Wardrop equilibrium in network games. At a Markovian network equilibrium, each player of the game solves a Markov decision process instead of a shortest path problem. We propose two novel extensions of Markovian network equilibrium by considering 1) variable demand, which offers the players a quitting option, and 2) multi-commodity flow, which allows players to have heterogeneous ending time. We further develop dynamic-programming-based iterative algorithms for the proposed equilibrium problems, together with their arithmetic complexity analysis. Finally, we illustrate our network equilibrium model via a multi-commodity ride-sharing example, and compare the computational efficiency of our algorithms against state-of-the-art optimization software Mosek over extensive numerical experiments.

\end{abstract}

\end{frontmatter}

\section{Introduction}
Network equilibrium problems arise in a variety of applications, such as resource allocation and routing in communication or transportation networks \cite{rockafellar1984network,bertsekas1998network,xiao2004simultaneous,burger2014duality}. Among the most well-studied examples is the Wardrop equilibrium model in routing games \cite{beckmann1956studies,gartner1980optimal1,gartner1980optimal2,correa2010wardrop,patriksson2015traffic}. In this model, users in a transportation network are assumed to choose routes with cost that they perceive as the lowest, \ie, each user solves a shortest path problem, under the prevailing traffic conditions \cite{correa2010wardrop}. With this assumption, the resulting equilibra are characterized by the Wardrop equilibrium principle: the cost of all the routes actually used are equal, and less than those which would be experienced by a single user on any unused route \cite{wardrop1952correspondence}.

To ensure their practical relevance, it is often necessary to incorporating stochasticity into the network equilibrium problems. For example, the stochastic user equilibrium (SUE) model \cite{fisk1980some,sheffi1982algorithm,liu2009method} considers independent stochastic error on the route cost perceived by the users, leading to user distribution based on the logit \cite{dial1971probabilistic} or probit model \cite{daganzo1977stochastic}; see \cite[Sec. 2.8.1]{patriksson2015traffic} and \cite{cominetti2012wardrop} for a detailed discussion. Unfortunately, the SUE model presents several drawbacks: it requires computationally expensive route enumeration, and is not suited for problems with overlapping routes due to its assumption of independent route cost. 

To address these drawbacks, different network models consider different type of stochasticity. In particular, \cite{baillon2008markovian,ahipacsaouglu2019distributionally} introduced a Markovian network equilibrium model where users are assumed to choose, instead of routes, sequences of actions with accumulated cost that they perceive as the lowest. Each action is accompanied by a deterministic outcome and a stochastic cost. For example, each vehicle in a transportation network is assumed to choose a sequence of arcs, where each arc leads to deterministic transition to the next node in the network and a stochastic amount of travel time \cite{baillon2008markovian}.

On the other hand, \cite{calderone2017markov,calderone2017infinite} proposed a different stochastic network equilibrium model. Unlike the one in \cite{baillon2008markovian}, each action is accompanied by a stochastic outcome and a deterministic cost. For example, an aircraft flying in stormy weather is assumed to choose a sequence of waypoints to fly towards, where each choice costs a deterministic amount of fuel usage and is accompanied by a stochastic change in the weather condition \cite{nilim2005robust}. As a result, instead of a shortest path problem, each user solves a Markov decision process (MDP) \cite{puterman2014markov,bertsekas1996neuro}, where the cost of different actions is determined by the prevailing choices of all users. This model has found a variety of applications in modern transportation including ridesharing and parking \cite{calderone2017models}. 

Although the results in \cite{calderone2017markov,calderone2017infinite} serves as a first step toward a more general class of stochastic dynamic network equilibrium model, it has the following limitations: 
a) it does not incorporate many important features of Wardrop equilibrium, such as variable demand and multi-commodity flow
and b) its solution method relies exclusively on off-the-shelf optimization software, which does not fully exploit the problem structure. We address these limitations by making the following contributions. 
\begin{enumerate}
    \item We develop novel extensions to the Markovian network equilibrium model by considering a) variable demand, which offers the users a quitting option, and b) multi-commodity flow, which allows users having heterogeneous ending time.
    \item We design novel dynamic-programming-based algorithms for Markovian network equilibrium problems with detailed arithematical complexity analysis. Our algorithms outperform state-of-the-art optimization software Mosek in extensive numerical experiments.  
\end{enumerate}

The rest of the paper is organized as follows. We first revisit some background on MDP in Section~\ref{section: preliminaries}, then present our variable demand and multi-commodity flow equilibrium models in Section~\ref{section: models}. Section~\ref{section: algorithms} focuses on developing efficient iterative algorithms for our equilibrium problems. Section~\ref{section: numerical examples} first illustrates the equilibrium models in Section~\ref{section: models} via a multi-commodity ride-sharing example, then compares the algorithms in Section~\ref{section: algorithms} against commercial software Mosek. Finally, we conclude with discussions and comments on the future directions of research in Section \ref{section: conclusion}. 

Throughout the paper we will use the following notation: \(\mathbb{R}\) denotes the set of real numbers, \(\mathbb{R}_+\) denotes the set of nonnegative real numbers, and \(\mathbb{N}\) denotes the set of positive integers; \([N]\) denotes the set \(\{1, 2, \ldots, N\}\) for integer \(N\); \(a_{ijk}\) denotes the \(ijk\)--th component of the three-dimensional tensor \(a\in\mathbb{R}^{n_1\times n_2\times n_3}\), and analogously, \(a_{ij}\) for the two-dimensional case. Given \(b_1, \ldots, b_N\in\mathbb{R}\), we say \((b^\star, i^\star)=\underset{i\in[N]}{\min}\,\,b_i,\) if \(
 b^\star= \underset{i\in[N]}{\min}\,\,b_i\) and   \(i^\star\in\underset{i\in[N]}{\argmin}\,\,b_i\).

\section{Preliminaries and background}
\label{section: preliminaries}
\begin{figure}
    \centering
    \def\layersep{3cm}

\begin{tikzpicture}[shorten >=1pt,->,draw=black!100, node distance=\layersep,scale=0.7]
    \tikzstyle{every pin edge}=[<-, shorten <=1pt]
    \tikzstyle{Rnode}=[rectangle,fill=black!25,minimum size=12pt,inner sep=0pt]
    \tikzstyle{Cnode}=[circle,fill=black!25,minimum size=12pt,inner sep=0pt]
    \tikzstyle{state}=[Cnode, fill=white!100,draw=black];
    \tikzstyle{action}=[Rnode, fill=white!100,draw=black];
    \tikzstyle{annot} = [text width=4em, text centered]

\foreach \name / \y in {1,2,3,5}
    \node[state] (I-\name) at (0,-\y) {};

\node at (0, -4) { \(\vdots\)};
\node at (\layersep, -3.5) { \(\vdots\)};
\node at (2*\layersep, -4) { \(\vdots\)};
\node at (-1, -3) { {\dots}};
\node at (7, -3) { {\dots}};

\foreach \name / \y in {1,2,3,5}
    \node[state] (O-\name) at (2*\layersep,-\y cm) {};

\foreach \name / \y in {1,2,4}
    \path[yshift=-0.5cm]
        node[action] (H-\name) at (\layersep,-\y cm) {};

\foreach \source in {1,2,3,5}
     \foreach \dest in {1,2,4}
         \path (I-\source) edge (H-\dest);
\foreach \source in {1,2,4}
    \foreach \dest in {1,2,3,5}
       \path (H-\source) edge (O-\dest);

\node[annot,above of=H-1, node distance=1.4cm] (hl) {\small actions};
\node[annot,left of=hl, node distance=2.2cm] {\small states};
\node[annot,right of=hl, node distance=2.2cm] {\small states};

\draw[draw=black, dashed] (-0.5,-0.5) rectangle ++(4,-5.5);
\node at (1.5, -5.6) {\small layer \(t\)};

\end{tikzpicture}
     
     \caption{A Markovian network}
     \label{fig: network of MDP}
 \end{figure}
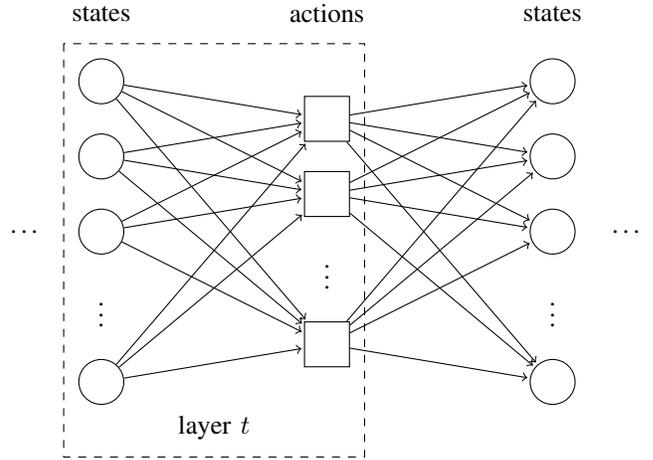

A \(T\)-horizon MDP is defined by a set of states \([S]\), a set of actions \([A]\), a cost tensor \(c\in\mathbb{R}^{T\times S\times A}\), and a transition probablity tensor \(P\in[0, 1]^{S\times A\times S}\), where \(T, S, A\in\mathbb{N}\) denote the number of time steps, states, and actions, respectively. Further, \(c_{tsa}\in\mathbb{R}\) denotes the cost of choosing action \(a\in[A]\) in state \(s\in[S]\) at time \(t\in[T]\), and \(P_{sas'}\in[0, 1]\) denotes the probability of transition from state \(s\in[S]\) to \(s'\in[S]\) when choosing action \(a\in[A]\). In order to find the optimal sequence of action that minimizes the expected accumulated cost, one can solve either one of the two following linear programs\footnote{Compared with the formulation in \cite{puterman2014markov}, the linear program here also allows \(p_{ts}>0\) when \(t>1\).}:
\begin{equation}
\begin{array}{ll}
\underset{y}{\mbox{min}} &  \sum\limits_{t, s, a} c_{tsa}y_{tsa}\\
\mbox{s.t.} &\sum\limits_{a} y_{1sa}=p_{1s},\\
& \sum\limits_{a}y_{t+1,sa}=p_{t+1,s}+\sum\limits_{s',a}P_{s'as}y_{ts'a},\enskip t\in[T-1],\\
&0\leq y_{tsa}, \enskip  \forall t\in[T], s\in[S], a\in[A].
\end{array}
\label{opt: linear optimal distribution}
\end{equation}
\begin{equation}
\begin{array}{ll}
\underset{v}{\mbox{max}} &  \sum\limits_{t,s}p_{ts} v_{ts}\\
\mbox{s.t.} & v_{Ts}\leq c_{Tsa},\\
&v_{ts}\leq c_{tsa} +\sum\limits_{s'}P_{sas'}v_{t+1, s'},\enskip t\in[T-1],\\
&\forall s\in[S], a\in[A].
\end{array}
\label{opt: linear optimal differential}
\end{equation}
where \(p\in\mathbb{R}_+^{T\times S}\) is such that \(p_{1s}>0\) for some \(s\in[S]\). If \(\sum_{s\in[S]}p_{1s}=1\) and \(p_{ts}=0\) for all \(1\leq t\leq T\) and \(s\in[S]\), then \(p_{1s}\) represents the probability of starting the MDP in state \(s\). Variable \(y_{tsa}\) in optimization \eqref{opt: linear optimal distribution} represents the probability of choosing action \(a\) in state \(s\) at time \(t\), and variable \(v_{ts}\) in optimization \eqref{opt: linear optimal differential} represents the expected accumulated cost between time \(t\) and time \(T\) starting from state \(s\). In general, optimization \eqref{opt: linear optimal distribution} and \eqref{opt: linear optimal differential} can be interpreted as the linear optimal distribution \cite[Sec.7A]{rockafellar1984network} and differential problem \cite[Sec.7E]{rockafellar1984network} defined on a \(T\)-layered \emph{Markovian network} (see Fig.~\ref{fig: network of MDP}), where \(p_{ts}\) represents the divergence on state node \(s\) in the \(t\)-th layer, \(y_{tsa}\) represents the flow from state \(s\) to action \(a\) in the \(t\)-th layer, and \(v_{ts}\) represents the potential on state \(s\) in the \(t\)-th layer. 

The following lemma shows that solutions of optimizations \eqref{opt: linear optimal distribution} and \eqref{opt: linear optimal differential} satisfy the dynamic programming principle. 

\begin{lemma}[\cite{puterman2014markov}]
\label{lem: MDP}
Suppose \(y\) solves \eqref{opt: linear optimal distribution}, and \(v\) solves \eqref{opt: linear optimal differential}. If \(y_{tsa}>0\) for any \(t\in[T], s\in[S], a\in[A]\), then
\begin{equation*}
\begin{aligned}
    &\textstyle (v_{Ts}, a)=\underset{a'\in[A]}{\min} \enskip c_{Tsa'}, \\
    &\textstyle (v_{ts}, a)=\underset{a'\in[A]}{\min}\enskip c_{tsa'}+\sum\limits_{s'}P_{sa's'}v_{t+1, s'},
    \end{aligned}
\end{equation*}
for all \(t\in[T-1]\).
\end{lemma}
Perhaps the most efficient solution algorithm for problem \eqref{opt: linear  optimal distribution} and \eqref{opt: linear optimal differential} is \emph{dynamic programming}, given by the following Algorithm~\ref{alg: Bellman} and Algorithm~\ref{alg: Kolmogorov}.  
\begin{algorithm}[!h]
\caption{Backward induction}
\begin{algorithmic}[1]
\Require \(P\), \(c\), \(T\).
\Ensure \(v, \pi\).
\State{Let \((v_{Ts}, \pi_{Ts})= \underset{a\in[A]}{\min} \,\,c_{Tsa},\enskip \forall s\in[S]\).}\label{alg: Bellman initial}
\For{\(t=T-1, T-2, \ldots, 1\)}
\State{\(\displaystyle (v_{ts}, \pi_{ts}) = \underset{a\in[A]}{\min}\big( c_{tsa}+\textstyle\sum\limits_{s'}P_{sas'}v_{t+1, s'}\big), \enskip \forall s\in[S]\)}\label{alg: dp induction step}
\EndFor
\end{algorithmic}
\label{alg: Bellman}
\end{algorithm}
\begin{algorithm}[!h]
\caption{Forward induction}
\begin{algorithmic}[1]
\Require \(\pi\), \(p\), \(P\), \(T\).
\Ensure \(y\).
\State{Initialize \(y=0\), let \(y_{1s\pi_{1s}}\gets
p_{1s}\) for all \(s\in[S]\).}

\For{\(t=1, 2, \ldots, T-1\)}
\State{\(y_{t+1, s\pi_{t+1, s}}\gets
p_{t+1, s}+\sum\limits_{j}P_{j\pi_{tj}s}y_{tj\pi_{tj}}, \enskip \forall s\in[S]\)}
\EndFor
\end{algorithmic}
\label{alg: Kolmogorov}
\end{algorithm}

Let \((v, \pi)\) be the output of Algorithm~\ref{alg: Bellman} with input \((P, c, T)\), and \(y\) be the output of Algorithm~\ref{alg: Kolmogorov} with input \((\pi, p, P, T)\), then one can easily show that such solution pair \((y, v)\) directly satisfies the Karush-Kuhn-Tucker (KKT) conditions \cite[Thm. 28.3]{rockafellar1970convex} of \eqref{opt: linear optimal distribution} and \eqref{opt: linear optimal differential}, hence it is an optimal primal-dual solution pair. If we define the \emph{sparsity level} of an MDP as follows
\begin{equation}\label{eqn: sparsity}
    \sigma=\max\{N_1, N_2\}/S\in [1/S, 1],
\end{equation}
where 
\(N_1=\max_{s, a}\big|\{s'|P_{sas'}>0\}\big|\) and \(N_2 =\max_{s', a}\big|\{s|P_{sas'}>0\}\big|\),
then \(\sigma S\) measures the maximum number of states connected by the transition kernel \(P\). Further, it is straightforward to check that Algorithm~\ref{alg: Bellman} costs \(O(\sigma TS^2A)\) arithmetic operations, and Algorithm~\ref{alg: Kolmogorov} costs \(O(\sigma TS^2)\) arithmetic operations. In addition, Algorithm~\ref{alg: Bellman} and Algorithm~\ref{alg: Kolmogorov} can be implemented as convolutional neural networks that allows efficient parallel computation \cite{tamar2016value}.

\section{Markovian network equilibrium}
\label{section: models}
By combining MDP together with classical routing games, Calderone and Sastry \cite{calderone2017markov} proposed MDP routing games where a fixed amount of players with the same planning horizon choose sequences of actions that they perceive as achieving the lowest expected accumulated cost under the prevailing choices of other players. Such games are similar to routing games where a fixed amount of players with the same destination choose routes that they perceive as the shortest under the prevailing choices of other players. 

In this section, we introduce two generalizations to MDP routing games that allow the amount of players to vary and the planning horizon to differ. We will also show that, under mild assumptions, the equilibra of such games can be computed efficiently using convex optimization. 

\subsection{Variable demand}
\label{subsection: variable demand}
One limitation of the MDP routing games in \cite{calderone2017markov} is the assumption that total amount of players is fixed. However, an important feature in network games is to allow the total amount of players to vary, or equivalently, to provide the players with a quitting action \cite[Sec. 2.1.2]{patriksson2015traffic}. Aiming to address this limitation, we propose the following variable demand MDP routing games.

\begin{game}\label{game: variable demand}
At each time \(t\in [T]\), \(p_{ts}\) new players start the game from state \(s\in [S]\). Among these \(p_{ts}\) players, each one can choose to
\begin{enumerate}
    \item quit the game immediately at the cost of \(\psi_{ts}(z_{ts})\),
    \item take action \(a \in[A]\) at the cost of \(\phi_{tsa}(y_{tsa})\) and reach state \(s'\in[S]\) with probability \(P_{sas'}\) at time \(t+1\), then repeat such process till \(t=T\), when the player ends the game after choosing the last action,
\end{enumerate}
where \(z_{ts}\) and \(y_{tsa}\) denote the total amount of players choosing to quit the game in state \(s\) at time \(t\), and, respectively, taking action \(a\) in state \(s\) at time \(t\).
\end{game}

\begin{remark}
Game~\ref{game: variable demand} is a special case of mean field games over graphs \cite{gomes2009discrete,gomes2010discrete,gueant2011infinity,gueant2015existence,tanaka2020linearly}. The interactions among different players is mediated by a mean field, described by function \(\phi_{tsa}\) and function \(\psi_{ts}\) for all \(t\in[T], s\in[S], a\in[A]\).
\end{remark}

Intuitively, one can interpret Game~\ref{game: variable demand} as a competitive market model. The supply side corresponds to the stochastic environment, providing the option of playing or quitting the game. The demand side corresponds to the amount of players that decided to play the game, which changes with the expected accumulated cost of the playing option according to curve \(\psi_{ts}\) for all \(t\in[T]\) and \(s\in[S]\).  

\begin{remark}
If the quitting option is not available, then Game~\ref{game: variable demand} reduces to an MDP routing game with fixed demand, introduced in \cite{calderone2017markov}. On the other hand, if the transition is Game~\ref{game: variable demand} is deterministic, \ie, for each \(s\in[S]\) and \(a\in[A]\), there exists \(s'\in[S]\) such that \(P_{sas'}=1\), then  Game~\ref{game: variable demand} reduces to a classical single-commodity routing game, with  and a variable demand \cite[Sec. 2.2.3]{patriksson2015traffic}. Particularly, each player solves an MDP with deterministic transition, which is equivalent to a shortest path problem.
\end{remark}

The Wardrop equilibrium principle is a key characterization of the equilibra of network games \cite{patriksson2015traffic,correa2010wardrop}. The principle states that, at equilibra, only the strategies with the lowest cost are actually used. Does this principle apply to Game~\ref{game: variable demand}? As we show in the following, the answer is affirmative. 

First, we make the following assumptions on Game~\ref{game: variable demand}.

\begin{assumption}\label{asp: variable demand} We assume that \(p\in\mathbb{R}_+^{T\times S}\), \(P\in[0, 1]^{S\times A\times S}\) and \(\sum_{s'}P_{sas'}=1\) for all \(s\in[S]\), \(a\in[A]\).
Further, the function \(\phi_{tsa}:[0, \rho]\to \mathbb{R}\) and function \(\psi_{ts}:[0, \rho]\to\mathbb{R}\) are continuous and strictly increasing over their respective domains, where \(\rho= \sum_{t, s}p_{ts}\). 
\end{assumption}

With these assumptions, we now introduce the following pair of primal-dual optimization problems associated with Game~\ref{game: variable demand}. 
\begin{equation}
\begin{array}{ll}
\underset{y, z}{\mbox{min}} & \sum\limits_{t,s,a} {\displaystyle\int_0^{y_{tsa}}\phi_{tsa}(\alpha)d\alpha}+\sum\limits_{t,s}\displaystyle\int_{0}^{z_{ts}}\psi_{ts}(\alpha)d\alpha\\
\mbox{s.t.} &\sum\limits_{a} y_{1sa}=p_{1s}-z_{1s},\\
& \sum\limits_{a}y_{t+1,sa}=p_{t+1, s}-z_{t+1, s}+\sum\limits_{s',a}P_{s'as}y_{ts'a},\\
&t\in [T-1],\\
&0\leq y_{tsa}, 0\leq z_{ts}\leq p_{ts}, \, \forall t\in[T], s\in[S], a\in[A]. 
\end{array}
\label{opt: optimal flow}
\end{equation}
\begin{equation}
\begin{array}{ll}
\underset{u, v, w, \lambda}{\mbox{max}} &
\sum\limits_{t,s}p_{ts} (v_{ts}-\lambda_{ts})-\sum\limits_{t,s,a} \displaystyle\int_{\phi_{tsa}(0)}^{u_{tsa}}\phi_{tsa}^{-1}(\alpha)d\alpha\\
&-\sum\limits_{t,s}\displaystyle\int_{\psi_{ts}(0)}^{w_{ts}}\psi^{-1}_{ts}(\alpha)d\alpha\\

\mbox{s.t.}  
&v_{Ts}\leq u_{Tsa},\\
&v_{ts}\leq u_{tsa} +\sum\limits_{s'}P_{sas'}v_{t+1, s'},\enskip  t\in [T-1],\\
&v_{ts}\leq w_{ts}+\lambda_{ts},\enskip 0\leq \lambda_{ts},\\
&\forall t\in[T],s\in[S], a\in[A].
\end{array}
\label{opt: optimal potential}
\end{equation}
In particular, the constraint \(0\leq z_{ts}\leq p_{ts}\) allows the number of players choosing to quit the game in state \(s\) at time \(t\) to vary winthin interval \([0, p_{ts}]\). If  variable \(z_{ts}\) is zero and function \(\phi_{tsa}\) is constant-valued for all \(t\in[T], s\in[S], a\in[A]\), \ie, the quitting option is removed and the cost of each action does not depend on \(y\) in in Game~\ref{game: variable demand}, then one can verify that optimization \eqref{opt: optimal flow} will reduce to \eqref{opt: linear optimal distribution} and optimization \eqref{opt: optimal potential} will reduce to \eqref{opt: linear optimal differential}.

The following theorem shows that, under Assumption~\ref{asp: variable demand}, the solution to the optimizations in \eqref{opt: optimal flow} and \eqref{opt: optimal potential} satisfy an equilibrium condition of Game~\ref{game: variable demand}. Similar to the Wardrop equilibrium principle, this equilibrium condition impliesthat  no individual player can benefit from unilaterally switching its actions.

\begin{theorem}\label{thm: variable demand} Suppose Assumption~\ref{asp: variable demand} holds, \((y, z)\) solves  \eqref{opt: optimal flow}, and \((u, v, w, \lambda)\) solves \eqref{opt: optimal potential}, then for any \(p_{ts}>0\),
\begin{equation}\label{eqn: Wardrop quitting}
\begin{aligned}
    &\text{if } z_{ts}=0,\enskip \text{then }v_{ts}\leq \psi_{ts}(p_{ts}),\\
    &\text{if } 0<z_{ts}< p_{ts}, \enskip \text{then }v_{ts}= \psi_{ts}(z_{ts}),\\
    &\text{if } z_{ts}= p_{ts}, \enskip \text{then }v_{ts}\geq \psi_{ts}(0).
    \end{aligned}
\end{equation}
Further, if \(y_{tsa}>0\), then
\begin{equation}\label{eqn: Wardrop playing}
\begin{aligned}
    &\textstyle (v_{Ts}, a) = \underset{a'\in[A]}{\min} \enskip \phi_{Tsa'}(y_{Tsa'}),\\
    &\textstyle (v_{ts}, a) = \underset{a'\in[A]}{\min}\enskip \phi_{tsa'}(y_{tsa'})+\sum\limits_{s'}P_{sa's'}v_{t+1,s'},
    \end{aligned}
\end{equation}
for all \(t\in[T-1]\).
\end{theorem}
\begin{pf*}{Proof}
See Appendix~\ref{app: thm variable demand}.
\end{pf*}

Theorem~\ref{thm: variable demand} shows that an equilibrium of Game~\ref{game: variable demand} that satisfies the Wardrop equilibrium principle not only exists, but can be computed by solving optimization \eqref{opt: optimal flow} and \eqref{opt: optimal potential}. In particular, if action \(a\) is chosen in state \(s\) at time \(t\) by any player at equilibrium, \ie,  \(y_{tsa}>0\), then action \(a\) must be optimal in the sense of Algorithm~\ref{alg: Bellman}. On the other hand, equation \eqref{eqn: Wardrop quitting} says that if some players choose the quitting option in state \(s\) at time \(t\) at equilibrium, \ie, \(z_{ts}>0\), then the cost of playing is no more than quitting, \ie, \(v_{ts}\geq \psi_{ts}(p_{ts})\). Similarly, if some players choose to play, \ie, \(z_{ts}<p_{ts}\), then the cost of playing is no more than quitting, \ie, \(v_{ts}\leq \psi_{ts}(p_{ts})\). Therefore, Theorem~\ref{thm: variable demand} indeed describes a Wardrop equilibrium where no individual player can benefit from unilaterally switching to alternative actions.

\subsection{Multicommodity flow}
\label{subsection: models: multicommodity flow}
Another limitation of the MDP routing games in \cite{calderone2017markov} is that all players are assumed to end their game at the same time, which is analogous to the single commodity routing game where all vehicles have the same destination. Aiming to address this limitation, we propose the following multi-commodity MDP routing game, where players can have heterogeneous ending time, denote by \(\mathbb{T}\). We assume, without loss of generality, that \(\mathbb{T}\subset[T]\) and \(T\in\mathbb{T}\). 

\begin{game}\label{game: multicommodity}
At each time \(t\in[T]\), \(p_{ts}^\tau\) new players who have a common ending time \(\tau\in\mathbb{T}\) with \(\tau\geq t\), start the game from state \(s\). Each of these players can choose the action \(a\) at the cost of \(\phi_{tsa}(\sum_{\tau, \tau\geq t}y_{tsa}^\tau)\) and reach state \(s'\) with probability \(P_{sas'}\) at time \(t+1\), then repeat this process till \(t=\tau\), when the player ends the game after choosing the last action. Here \(y_{tsa}^\tau\) denotes the total amount of players who plan to end the game at time \(\tau\) and choose action \(a\) in state \(s\) at time \(t\).
\end{game}

\begin{remark}
If \(\mathbb{T}= \{T\}\), then Game~\ref{game: multicommodity} reduces to a MDP routing game introduced in \cite{calderone2017markov}. On the other hand, if the transition in Game~\ref{game: multicommodity} is deterministic, \ie, for each \(s\in[S]\) and \(a\in[A]\), there exists \(s'\in[S]\) such that \(P_{sas'}=1\), then Game~\ref{game: multicommodity} reduces to the traditional multi-commodity routing game with a fixed demand \cite[Sec. 2.1.1]{patriksson2015traffic}. Particularly, the state where a player start and end the game form a origin-destination pair, which is jointly determined by the starting state and the deterministic transition.
\end{remark}

Similar to Game~\ref{game: variable demand}, the equilibrium of Game~\ref{game: multicommodity} can also be computed by solving convex optimization problems, as we show in the following.

First, we make the following assumptions on Game~\ref{game: multicommodity}.
\begin{assumption}\label{asp: multicommodity} We assume \(T\in\mathbb{T}\subseteq[T]\), \(p_{ts}^\tau\in\mathbb{R}_+\) for all \(\tau\in\mathbb{T}\), \(t\geq \tau\) and \(s\in[S]\); \(P\in\mathbb{R}_+^{S\times A\times S}\) and \(\sum_{s'}P_{sas'}=1\) for all \(s\in[S]\), \(a\in[A]\). Further, the
function \(\phi_{tsa}:[0, \rho]\to \mathbb{R}\) is continuous and strictly increasing, where \(\rho = \sum_\tau\sum_{t\leq \tau, s}p_{ts}^\tau\).
\end{assumption}

With these assumptions, we now introduce the following pair of primal-dual optimization problems associated with Game~\ref{game: multicommodity}. Notice that if \(\mathbb{T}=\{T\}\), then they reduce to optimization \eqref{opt: linear optimal distribution} and \eqref{opt: linear optimal differential}, respectively. 
\begin{equation}
    \begin{array}{ll}
        \underset{\{y^\tau\}_{\tau\in\mathbb{T}}}{\mbox{min}} & \sum\limits_{t,s,a} \displaystyle\int_0^{\sum\limits_{\tau, \tau\geq t}y^{\tau}_{tsa}} \phi_{tsa}(\alpha)d\alpha \\
         \mbox{s.t.} & \sum\limits_{a} y^{\tau}_{1sa}=p^{\tau}_{1s}\\ &\sum\limits_{a}y^{\tau}_{t+1, sa}=p^{\tau}_{t+1,s}+\sum\limits_{s',a}P_{s'as}y^{\tau}_{ts'a},\, t\in[\tau-1],\\
         &0\leq y^{\tau}_{tsa}, \enskip
         \forall \tau\in \mathbb{T}, t\in[\tau], s\in[S], a\in[A]
    \end{array}
    \label{opt: multi optimal flow}
\end{equation}
\begin{equation}
    \begin{array}{ll}
        \underset{u, \{v^\tau\}_{\tau\in\mathbb{T}}}{\mbox{max}} &  \sum\limits_{t,s}\,\,\sum\limits_{\tau, \tau\geq t} p^{\tau}_{ts}v^{\tau}_{ts}-\sum\limits_{t,s,a} \displaystyle\int_{\phi_{tsa}(0) }^{u_{tsa}} \phi_{tsa}^{-1}(\alpha)d\alpha \\
         \mbox{s.t.} & v^{\tau}_{\tau s}\leq u_{\tau sa},\\
         &v^{\tau}_{ts}\leq u_{tsa}+\sum\limits_{s'}P_{sas'}v^{\tau}_{t+1, s'}, \enskip t\in[\tau-1],\\
         &\forall \tau\in\mathbb{T}, s\in[S], a\in[A]
    \end{array}
    \label{opt: multi optimal potential}
\end{equation}

The following theorem shows that, under Assumption~\ref{asp: multicommodity}, the solutions to optimization problems \eqref{opt: multi optimal flow} and \eqref{opt: multi optimal potential} satisfy the equilibrium condition of Game~\ref{game: multicommodity}. Similar to the Wardrop equilibrium principle, this equilibrium condition implies that no individual player can benefit from unilaterally switching actions. 

\begin{theorem}\label{thm: multicommodity} Suppose Assumption~\ref{asp: multicommodity} holds, \(y\) solves \eqref{opt: multi optimal flow}, and \((u, v)\) solves \eqref{opt: multi optimal potential}. If \(y^{\tau}_{tsa}>0\) for any \( \tau\in\mathbb{T}, s\in[S], a\in[A]\), then
\begin{equation}
\begin{aligned}
    &\textstyle (v^\tau_{\tau s}, a)=\underset{a'\in[A]}{\min} \enskip \phi_{\tau sa'}\big(\sum\limits_{\tau, \tau\geq t}y^{\tau}_{\tau sa'}\big), \\
    &\textstyle (v^\tau_{ts}, a)=\underset{a'\in[A]}{\min}\enskip \phi_{tsa'}\big(\sum\limits_{\tau, \tau\geq t}y^{\tau}_{tsa'}\big)+\sum\limits_{s'}P_{sa's'}v^{\tau}_{t+1, s'},
    \end{aligned}\label{eqn: multi Wardrop}
\end{equation}
for all \(t\in[\tau-1]\).
\end{theorem}
\begin{pf*}{Proof}
See Appendix~\ref{app: thm multicommodity}.
\end{pf*}

Theorem~\ref{thm: multicommodity} shows that a Wardrop equilibrium of Game~\ref{game: multicommodity} not only exists, but can be found by solving optimization problems \eqref{opt: multi optimal flow} and \eqref{opt: multi optimal potential}. In particular, the equations in \eqref{eqn: multi Wardrop} characterize a multi-commodity flow Wardrop equilibrium in the sense that no individual player can benefit from using alternative actions before his/her ending time \(\tau\) for all \(\tau\in\mathbb{T}\).

\section{Efficient algorithms via linearization}
\label{section: algorithms}
In this section, we develop efficient iterative algorithms for the network equilibrium problems introduced in the previous section. In particular, we first prove that the linearized versions of problem \eqref{opt: optimal flow} and problem \eqref{opt: multi optimal flow} can both be solved in closed form via Algorithm~\ref{alg: Bellman} and Algorithm~\ref{alg: Kolmogorov}. This observation motivates efficient iterative algorithms that enjoy detailed arithematical complexity analysis. 

We will use the following notation to simply our later discussions. Given \(y, u\in\mathbb{R}^{T\times S\times A}\) and \(z, w\in\mathbb{R}^{T\times S}\), we let \(\phi(y), \phi^{-1}(u)\in\mathbb{R}^{T\times S\times A}\) and \(\psi(z), \psi^{-1}(w)\in\mathbb{R}^{T\times S}\) be such that
\begin{equation}\label{eqn: compact tensor}
\begin{aligned}
    &[\phi(y)]_{tsa}=\phi_{tsa}(y_{tsa}), \enskip [\phi^{-1}(u)]_{tsa}=\phi^{-1}_{tsa}(u_{tsa}),\\
    &[\psi(z)]_{ts}=\phi_{ts}(z_{ts}),\enskip  [\psi^{-1}(w)]_{ts}=\psi^{-1}_{ts}(w_{ts}),
\end{aligned}
\end{equation}
for all \(t\in[T], s\in[S], a\in[A]\). We also let \(\underline{u}, \overline{u}\in\mathbb{R}^{T\times S\times A}\) and \(\underline{w}, \overline{w}\in\mathbb{R}^{T\times S}\) be such that
\begin{equation}
    \begin{aligned}
        & \underline{u}_{tsa}=\phi_{tsa}(0), \enskip \overline{u}_{tsa}=\phi_{tsa}(\rho),\\
        & \underline{w}_{ts}=\psi_{ts}(0), \enskip \overline{w}_{ts}=\psi_{ts}(\rho),
    \end{aligned}
\end{equation}
for all \(t\in[T], s\in[S], a\in[A]\).

\subsection{Linearization and dynamic programming}

If we approximate the objective function in \eqref{opt: optimal flow} using its linearization at \(u\in\mathbb{R}^{T\times S\times A}\) and  \(w\in\mathbb{R}^{T\times S}\), we obtain the following
\begin{equation}
\begin{array}{lll}
-g(u, w)=&\underset{y, z}{\mbox{min}} & \sum\limits_{t,s,a} u_{tsa}y_{tsa}+\sum\limits_{t,s}w_{ts}z_{ts}\\
&\mbox{s.t.} & \text{constraints in problem \eqref{opt: optimal flow}. }
\end{array}
\label{opt: linear optimal flow}
\end{equation}

Observe that the above optimization is a modification to problem \eqref{opt: linear optimal distribution}, by including an additional variable \(z\). This suggest that \eqref{opt: linear optimal flow} may also be solved using dynamic programming, which is confirmed by the following lemma.

\begin{lemma}\label{lem: variable demand}
Suppose Assumption~\ref{asp: variable demand} holds. Let \((\hat{v}, \hat{\pi})\) be the output of Algorithm~\ref{alg: Bellman} with input \((P, u, T)\), and
\[\hat{z}_{ts}=\begin{cases}
p_{ts}, & \hat{v}_{ts}>w_{ts}\\
0, & \hat{v}_{ts}\leq w_{ts}
\end{cases}\enskip \forall t\in[T], s\in[S].\]
In addition, let \(\hat{y}\) be the output of Algorithm~\ref{alg: Kolmogorov} with input \((\hat{\pi}, p-\hat{z}, P, T)\). Then 
\[-g(u, w)=\textstyle\sum\limits_{t, s, a}u_{tsa}\hat{y}_{tsa}+\sum\limits_{t, s}w_{ts}\hat{z}_{ts}.\]
Further, for any \(u'\in\mathbb{R}^{T\times S\times A}\) and \(w'\in\mathbb{R}^{T\times S}\),
\begin{equation*}
\begin{aligned}
    &g(u', w')-g(u, w)\\
    \geq &\textstyle \sum\limits_{t, s, a}(u'_{tsa}-u_{tsa})(-\hat{y}_{tsa})+\sum\limits_{t, s}(w'_{ts}-w_{ts})(-\hat{z}_{ts}).
\end{aligned}
\end{equation*}

\end{lemma}
\begin{pf*}{Proof}
See Appendix~\ref{app: lem variable demand}.
\end{pf*}

Similarly, if we approximate the objective function in \eqref{opt: multi optimal flow} using a linear function, we obtain the following
\begin{equation}
    \begin{array}{lll}
       -h(u)=& \underset{y^\tau, \tau\in\mathbb{T}}{\mbox{min}} & \sum\limits_{t,s,a}\,\, \sum\limits_{\tau, \tau\geq t} u_{tsa}y_{tsa}^\tau \\
         &\mbox{s.t.} & \text{constraints in problem \eqref{opt: multi optimal flow}. }
    \end{array}
    \label{opt: linear multi optimal flow}
\end{equation}
where \(u\in\mathbb{R}^{T\times S\times A}\) is the approximation parameter, \(-h(u)\) is the optimal value of \eqref{opt: linear multi optimal flow}. The following lemma shows that optimization \eqref{opt: linear multi optimal flow} can be solved using Algorithm~\ref{alg: Bellman} and \ref{alg: Kolmogorov} as well.

\begin{lemma}\label{lem: multicommodity}
Suppose Assumption~\ref{asp: multicommodity} holds. Let \((\hat{v}^\tau, \hat{\pi})\) be the output of Algorithm~\ref{alg: Bellman} with input \((P, u, \tau)\), \(\hat{y}^\tau\) be the output of Algorithm~\ref{alg: Kolmogorov} with input \((\hat{\pi}^\tau, p^\tau, P, \tau)\). Then
\[-h(u)=\textstyle \sum\limits_{t,s,a}\sum\limits_{\tau, \tau\geq t} u_{tsa}\hat{y}^\tau_{tsa}.\]
Further, for any \(u'\in\mathbb{R}^{T\times S\times A}\),
\[h(u')-h(u)\geq \textstyle \sum\limits_{t,s,a}\sum\limits_{\tau, \tau\geq t} (u'_{tsa}-u_{tsa})(-\hat{y}^\tau_{tsa}).\]
\end{lemma}
\begin{pf*}{Proof}
See Appendix~\ref{app: lem multicommodity}.
\end{pf*}

\begin{remark}\label{rem: support function}
Function \(g:\mathbb{R}^{T\times S\times A}\times \mathbb{R}^{T\times S}\to\mathbb{R}\) in Lemma~\ref{lem: variable demand} is the support function of a polyhedron, which is closed and convex \cite[p.28]{rockafellar1970convex}. Further, Lemma~\ref{lem: variable demand} shows that the slope of a linear underestimator, or \emph{subgradient}, of function \(g(u, w)\) can be computed using Algorithm~\ref{alg: Bellman} and Algorithm~\ref{alg: Kolmogorov}. Similar observation is made in Lemma~\ref{lem: multicommodity} for function \(h:\mathbb{R}^{T\times S\times A}\to\mathbb{R}\). 
\end{remark}

\subsection{Iterative algorithms using linearization}
We now develop iterative algorithms for optimization problems in Section~\ref{section: models} using the results from the previous subsection. We will use the following notion of \(\epsilon\)-optimal solution.

\begin{definition}
Given a constrained optimization where an objective function is optimized subject to constraints, we say a solution is \(\epsilon\)-optimal \(\epsilon\in\mathbb{R}_+\) if it satisfies all the constraints and the objective function value evaluated at this solution is at most \(\epsilon\) away from the optimal value.
\end{definition}

We will also use the following additional assumptions on Game~\ref{game: variable demand} and, respectively, Game~\ref{game: multicommodity}. 

\begin{assumption}\label{asp: Lipschitz variable demand}
Function \(\phi_{tsa}:[0, \rho]\to\mathbb{R}\) and \(\psi_{ts}:[0, \rho]\to\mathbb{R}\) are \(L\)-Lipschitz continuous over their respective domains for all \(t\in[T]\), \(s\in[S]\) and \(a\in[A]\).
\end{assumption}

\begin{assumption}\label{asp: Lipschitz multicommodity}
Function \(\phi_{tsa}:[0, \rho]\to\mathbb{R}\) is \(L\)-Lipschitz continuous over its domain for all \(t\in[T]\), \(s\in[S]\) and \(a\in[A]\).
\end{assumption}

\begin{remark}\label{rem: Lipschitz}
Assumption~\ref{asp: Lipschitz variable demand} and Assumption~\ref{asp: Lipschitz multicommodity} are mild assumptions on the differentiability of the corresponding functions. For example, if \(\phi_{tsa}\) is continuously differentiable, then the mean value theorem states that for any \(\alpha_1, \alpha_2\in[0, \rho]\), there exists \(\alpha_3\in[0, \rho]\) such that
\[|\phi_{tsa}(\alpha_1)-\phi_{tsa}(\alpha_2)|\leq |\phi'_{tsa}(\alpha_3)|\cdot |\alpha_1-\alpha_2|.\]
where \(\phi'_{tsa}\) is the derivative of function \(\phi_{tsa}\). Hence Assumption~\ref{asp: Lipschitz multicommodity} is satisfied by choosing \[L\geq \max_{\alpha\in[0, \rho]} |\phi'_{tsa}(\alpha)|,\enskip \forall t\in[T], s\in[S], a\in[A].\]
which takes a bounded value since \(\phi'_{tsa}\) is continuous. However, the continuity of \(\phi'_{tsa}\) is not necessary. For example, if \(\phi_{tsa}\) is a piecewise linear function, \ie, a function that is affine over a collection of intervals, then it is still Lipschitz continuous even if its derivative is not continuous. 
\end{remark}

Based on Lemma~\ref{lem: variable demand} and Lemma~\ref{lem: multicommodity}, we propose to solve optimiztaion \eqref{opt: optimal flow} and \eqref{opt: multi optimal flow} using Frank-Wolfe method \cite{frank1956algorithm}, which repeatedly solve the linearized versions of \eqref{opt: optimal flow} and \eqref{opt: multi optimal flow}. We summarize the Frank-Wolf method for optimization \eqref{opt: optimal flow} and \eqref{opt: multi optimal flow} in Algorithm~\ref{alg: Frank-Wolfe} and, respectively, Algorithm~\ref{alg: multi Frank-Wolfe}. The following theorem provides the overall arithmetic complexity analysis of Algorithm~\ref{alg: Frank-Wolfe} and Algorithm~\ref{alg: multi Frank-Wolfe}.

\begin{algorithm}[!ht]
\caption{Frank-Wolfe method}
\label{alg: Frank-Wolfe}
\begin{algorithmic}[1]
\Require \(p, P, \phi, \psi, T, \{\alpha^k\}\), initial value for \(y, z\).
\For{\(k=1, 2, \ldots, K\)}
\State {\((\hat{v}, \hat{\pi})\gets\) Alg.~\ref{alg: Bellman}\((P, \phi(y), T)\).}
\State{\(\hat{z}_{ts}=\begin{cases}
    p_{ts}, & \hat{v}_{ts}>\psi_{ts}(z_{ts})\\
    0, & \hat{v}_{ts}\leq \psi_{ts}(z_{ts})
    \end{cases}, \enskip \forall t\in[T], s\in[S]\)}
\State{\(\hat{y}\gets\) Alg.~\ref{alg: Kolmogorov} \(( \hat{\pi}, p-\hat{z}, P, T)\).}
\State{\(y\gets y-\alpha^k(y- \hat{y})\)}
\State{\(z\gets z-\alpha^k(z- \hat{z})\)}
\EndFor
\end{algorithmic}
\end{algorithm}

\begin{algorithm}[!ht]
\caption{Multicommodity Frank-Wolfe method}
\label{alg: multi Frank-Wolfe}
\begin{algorithmic}[1]
\Require \(p, P, \phi, \mathbb{T}, \{\alpha^k\}\), initial value for \(y^\tau\) for all \(\tau\in\mathbb{T}\).
\For{\(k=1, 2, \ldots, K\)}
\State{\(\hat{y}_{tsa}=\sum_{\tau, \tau\geq t}y_{tsa}^{\tau}, \enskip \forall t\in[T], s\in[S], a\in[A]\)}
\State {\((\hat{v}^\tau, \hat{\pi}^\tau)\gets\)Alg.~\ref{alg: Bellman}\((P, \phi(\hat{y}), \tau), \enskip \forall \tau\in\mathbb{T}\)}
\State{\(\hat{y}^\tau\gets\) Alg.~\ref{alg: Kolmogorov}\((\hat{\pi}^\tau, p^\tau, P, \tau), \enskip \forall \tau\in\mathbb{T}\)}
\State{\(y^\tau\gets y^\tau-\alpha^k(y^\tau- \hat{y}^\tau), \enskip \forall \tau\in\mathbb{T}\)}
\EndFor
\end{algorithmic}
\end{algorithm}

The following theorem shows the convergence property of Algorithm~\ref{alg: Frank-Wolfe} and Algorithm~\ref{alg: multi Frank-Wolfe}.

\begin{theorem}\label{thm: FW convergence}
Let \(\sigma\) be given by \eqref{eqn: sparsity}. If Assumption~\ref{asp: variable demand} and \ref{asp: Lipschitz variable demand} hold, then Algorithm~\ref{alg: Frank-Wolfe} with \(\alpha^k=\frac{2}{k+1}\) gives an \(\epsilon\)-optimal solution to \eqref{opt: optimal flow} in \(O(\frac{\sigma TS^2A}{\epsilon})\) arithmetic operations. Similarly, if Assumption~\ref{asp: multicommodity} and \ref{asp: Lipschitz multicommodity} hold, then Algorithm~\ref{alg: multi Frank-Wolfe} with \(\alpha^k=\frac{2}{k+1}\) gives an \(\epsilon\)-optimal solution to \eqref{opt: multi optimal flow} in \(O(\frac{\sigma T^2S^2A}{\epsilon})\) arithmetic operations. 
\end{theorem}
\begin{pf*}{Proof}
See Appendix~\ref{app: thm FW convergence}.
\end{pf*}

Theorem~\ref{thm: variable demand} provides arithmetical complexity of Algorithm~\ref{alg: Frank-Wolfe} and \ref{alg: multi Frank-Wolfe}, which not only depends on the problem size (\ie, \(T, S, A\)),  but also the sparsity of the constraints (\ie, \(\sigma\)) in \eqref{opt: optimal flow} and \eqref{opt: multi optimal flow}.

What about the dual problems? Observe that the optimization in \eqref{opt: optimal potential} can be separated into two layers: an outer layer that optimizes over \((u, w)\), and an inner layer that optimizes over \((v, \lambda)\) for a given value of \((u, w)\); namely, \eqref{opt: optimal potential} is equivalent to the following
\begin{equation}
\begin{array}{ll}
\underset{u, w}{\mbox{max}} &
-g(u, w)-\sum\limits_{t,s,a} \displaystyle\int_{\phi_{tsa}(0)}^{u_{tsa}}\phi_{tsa}^{-1}(\alpha)d\alpha\\
&-\sum\limits_{t,s}\displaystyle\int_{\psi_{ts}(0)}^{w_{ts}}\psi^{-1}_{ts}(\alpha)d\alpha\\
\end{array}
\end{equation}
where
\begin{equation}
\begin{array}{lll}
-g(u, w)=&\underset{v, \lambda}{\mbox{max}} &
\sum\limits_{t,s}p_{ts} (v_{ts}-\lambda_{ts})\\
&\quad \mbox{s.t.}  & \text{constraints in \eqref{opt: optimal potential}.}
\end{array}
\label{opt: linear optimal potential}
\end{equation}
One can show that optimization \eqref{opt: linear optimal potential} is exactly the dual problem of \eqref{opt: linear optimal flow}. Since the constraint sets in \eqref{opt: linear optimal flow} and \eqref{opt: linear optimal potential} are both nonempty, the optimal value of \eqref{opt: linear optimal flow} and \eqref{opt: linear optimal potential} are the same \cite{morgenstern1953theory}.
In other words, \eqref{opt: optimal potential} can be equivalently written as follows 
\begin{equation}
\begin{array}{ll}
\underset{u, w}{\mbox{max}} &
-g(u, w)-\sum\limits_{t,s,a} \displaystyle\int_{\phi_{tsa}(0)}^{u_{tsa}}\phi_{tsa}^{-1}(\alpha)d\alpha\\
&-\sum\limits_{t,s}\displaystyle\int_{\psi_{ts}(0)}^{w_{ts}}\psi^{-1}_{ts}(\alpha)d\alpha\\
\mbox{s.t.}  
&  \text{\(-g(u, w)\) is the optimal value of \eqref{opt: linear optimal flow}. }
\end{array}
\label{opt: outer optimal potential}
\end{equation}

Using similar reasoning, we can rewrite \eqref{opt: multi optimal potential} as follows
\begin{equation}
    \begin{array}{ll}
        \underset{u}{\mbox{max}} &  -h(u)-\sum\limits_{t,s,a} \displaystyle\int_{\phi_{tsa}(0) }^{u_{tsa}} \phi_{tsa}^{-1}(\alpha)d\alpha \\
         \mbox{s.t.} & \text{ \(-h(u)\) is the optimal value of \eqref{opt: linear multi optimal flow}. }
    \end{array}
    \label{opt: outer multi optimal potential}
\end{equation}

We already discussed in Remark~\ref{rem: support function} how the subgradients of function \(g(u, w)\) and \(h(u)\) can be computed efficiently using Algorithm~\ref{alg: Bellman} and Algorithm~\ref{alg: Kolmogorov}. In addition, all the other terms in the objective functions of problem \eqref{opt: outer optimal potential} and \eqref{opt: outer multi optimal potential} are continuously differentiable. This suggest that \eqref{opt: outer optimal potential} and \eqref{opt: outer multi optimal potential} are suited for the projected subgradient method, which optimizes a non-smooth function by repeatedly computing its subgradients and projections onto its domain. We summarize the projected subgradient method applied to \eqref{opt: outer optimal potential} and \eqref{opt: outer multi optimal potential} in Algorithm~\ref{alg: subgradient} and, respectively, Algorithm~\ref{alg: multi subgradient}.

\begin{algorithm}
\caption{Subgradient method}
\begin{algorithmic}[1]
\Require \(P, p, \phi, \psi, T, \{\alpha^k\}\), initial value for \(u, w\).
\For{\(k=1, 2, \ldots, K\)}
\State {\((\hat{v}, \hat{\pi})\gets\) Alg.~\ref{alg: Bellman}\((P, u, T)\)}
\State{
\(\hat{z}_{ts}=\begin{cases}
    p_{ts}, & \hat{v}_{ts}>w_{ts}\\
    0, & \hat{v}_{ts}\leq w_{ts}
    \end{cases}, \enskip \forall t\in[T], s\in[S]\)} 
\State{\(\hat{y}\gets\) Alg.~\ref{alg: Kolmogorov} \(( \hat{\pi}, p-\hat{z}, P, T)\).}
\State{\(u\gets \min\{\overline{u},\max\{\underline{u},u+\alpha^k\left(\hat{y}- \phi^{-1}(u)\right)\}\}\)}
\State{\(w\gets \min\{\overline{w},\max\{\underline{w}, w+\alpha^k\left(\hat{z}-\psi^{-1}(w)\right)\}\}\)}
\EndFor
\end{algorithmic}
\label{alg: subgradient}
\end{algorithm}

\begin{algorithm}
\caption{Multi-commodity subgradient method}
\begin{algorithmic}[1]
\Require \(p, P, \phi, \mathbb{T}, \{\alpha^k\}\), initial value of \(u\).
\For{\(k=1, 2, \ldots, K\)}
\State {\((\hat{v}^\tau, \hat{\pi}^\tau)\gets\) Alg.~\ref{alg: Bellman}\((P, u, \tau), \enskip \forall\tau\in\mathbb{T}\)} 
\State{\(\hat{y}^\tau\gets\)Alg. \ref{alg: Kolmogorov}\((\hat{\pi}^\tau, p^\tau, P, \tau), \enskip \forall\tau\in\mathbb{T}\)}
\State{\(\hat{y}_{tsa}=\sum_{\tau, \tau\geq t}y_{tsa}^{\tau}, \enskip \forall t\in[T], s\in[S], a\in[A]\)}
\State{\(u\gets \min\{\overline{u},\max\{\underline{u}, u+\alpha^k\left(\hat{y}- \phi^{-1}(u)\right)\}\}\)}
\EndFor
\end{algorithmic}
\label{alg: multi subgradient}
\end{algorithm}

The following theorem shows the the convergence property of Algorithm~\ref{alg: subgradient} and Algorithm~\ref{alg: multi subgradient}.
\begin{theorem}\label{thm: subgrad convergence}
Let \(\sigma\) be given by \eqref{eqn: sparsity}. If Assumption~\ref{asp: variable demand} and Assumption~\ref{asp: Lipschitz variable demand} hold, then Algorithm~\ref{alg: multi subgradient} with \(\alpha^k=\frac{2L}{k+1}\) gives an \(\epsilon\)-optimal solution to \eqref{opt: outer optimal potential} using \(O(\frac{\sigma TS^2A}{\epsilon})\) arithmetic operations. Similarly, if Assumption~\ref{asp: multicommodity} and Assumption~\ref{asp: Lipschitz multicommodity} hold, then Algorithm~\ref{alg: multi subgradient} with \(\alpha^k\equiv\frac{2LT}{k+1}\) gives an \(\epsilon\)-optimal solution to \eqref{opt: outer multi optimal potential} in \(O(\frac{\sigma T^2S^2A}{\epsilon})\) arithmetic operations. 
\end{theorem}
\begin{pf*}{Proof}
See Appendix~\ref{app: subgrad convergence}
\end{pf*}

\section{Numerical examples}
\label{section: numerical examples}
 In this section, we first illustrate the equilibrium models in Section~\ref{section: models} using a ride-sharing example, then demonstrate the efficiency of algorithms in Section~\ref{section: algorithms} by comparing them against commercial software Mosek (\url{https://www.mosek.com}) over extensive numerical experiments.

\subsection{Multicommodity ride-sharing game}\label{subsection: ride-sharing}

\begin{figure}[ht]
    \centering
    \includegraphics[width=0.8\linewidth]{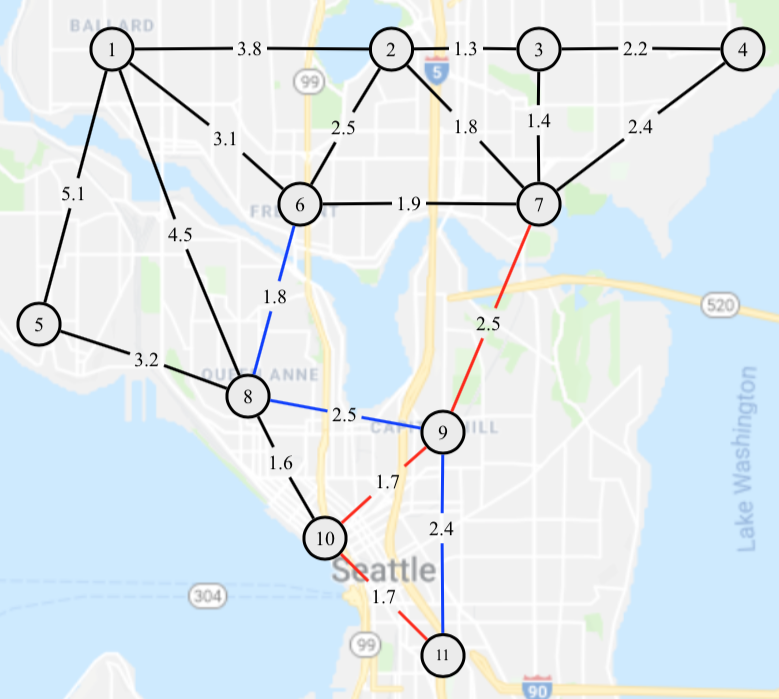}
    \caption{Seattle transportation network and candidate LRT routes: 7-9-10-11(red) and 6-8-9-11(blue).} 
    \label{fig: seattle}
\end{figure}

We consider the game played by ride-sharing drivers in Seattle, competing for customers. We first abstract the Seattle area as an undirected graph illustrated in Fig.~\ref{fig: seattle}, whose nodes denote various neighborhoods in Seattle, and edges denote available routes, labeled by its driving distance. We denote  the set of neighboring nodes of node \(s\) as \(\mathcal{N}_s\). We model the decision-making of an ride-sharing driver on a typical weekend night (7pm-1am) as an MDP defined as follows.

\begin{itemize}
    \item Time steps: \(t=1, 2, \ldots, 36\) denotes the (end of) 10-minute-intervals between 7pm and 1am.   
    \item States: \([S]\) correspond to different nodes in graph \(\mathcal{G}\). 
    \item Actions: in state \(s\), \(a_{s'}\) denotes picking up a waiting rider with destination \(s'\) for all \(s'\in\mathcal{N}_s\); \(a_{\text{wait}}\) denotes waiting for a future rider.
    \item Transition kernel: we assume \(P_{sas'}\) is given by
    \begin{equation*}
        P_{sas'}=\begin{cases}
    1, & \text{if \(a=a_{s'}\), \(s'\in \mathcal{N}_s\),}\\
    1/(|\mathcal{N}_s|+1), & \text{if \(a=a_{\text{wait}}\), \(s'\in \mathcal{N}_s\cup\{s\}\).}
    \end{cases}
    \end{equation*}
    All other entries of \(P_{sas'}\) are zero. Here we use an uniform distribution over neighboring states to describe the uncertain destinations of future riders\footnote{Such distribution can be approximated more accurately using historical data in practical applications.}.
    
    \item Cost: due to the competition among drivers, we assume the profit for picking up a rider decreases with the amount of drivers making the same offer, namely
    \begin{equation}\label{eqn: profit}
        f_{tss'}=\alpha+\beta   \big(1-\frac{y_{tsa_{s'}}}{\gamma_{tss'}}\big)\text{dist}_{ss'},
    \end{equation}
    for all \(t\in[T], s\in[S]\), where \(\alpha\) and \(\beta\) is the baseline profit and, respectively, nominal profit per mile. We let \(\text{dist}_{ss'}\) denotes distance(miles) between \(s\) and \(s'\), \(\gamma_{tss'}\) denotes the rider demand from \(s\) to \(s'\) at time \(t\), and finally \(y_{tsa_{s'}}\) denotes the amount of drivers choosing action \(a_{s'}\) in state \(s\) at time \(t\). The cost of action \(a\) in state \(s\) is a function of \(y_{tsa}\) defined as follows
    \begin{equation*}
        \phi_{tsa}(y_{tsa})=\begin{cases}
    -f_{tss'}, & \text{if \(a=a_{s'},s'\in\mathcal{N}_s\).}\\
    -\sum\limits_{s'\in\mathcal{N}_s}P_{sas'}f_{tss'}, & \text{if \(a=a_{\text{wait}}\).}
    \end{cases}
    \end{equation*}
    \item Planning time windows: We assume that 10 drivers start working from each state every 10 minutes between 7pm and 9pm. Once started, each driver is assumed to only work for 4 consecutive hours to avoid driver fatigue, \ie , \(p^{(t+24)}_{ts}=10\) for all \(s\in[S]\) and \(t\in[12]\). 
\end{itemize}
Notice that the function \(\phi_{tsa}\) defined above is linear with slope \(\alpha\), hence Assumption~\ref{asp: multicommodity} is satisfied with \(L=\alpha\). In general, as long as \(\phi_{tsa}\) is modeled or approximated as a continuously differentiable function, Assumption~\ref{asp: Lipschitz multicommodity} is always satisfied, as we discussed in Remark~\ref{rem: Lipschitz}. We also assume that each driver can travel between neighboring nodes within one time step in this simplified transportation network. In practice, such assumption can be ensured by adding more nodes to the network using a finer discretization of the interested area.

\begin{figure}[ht]
    \centering
    \includegraphics[width=\linewidth]{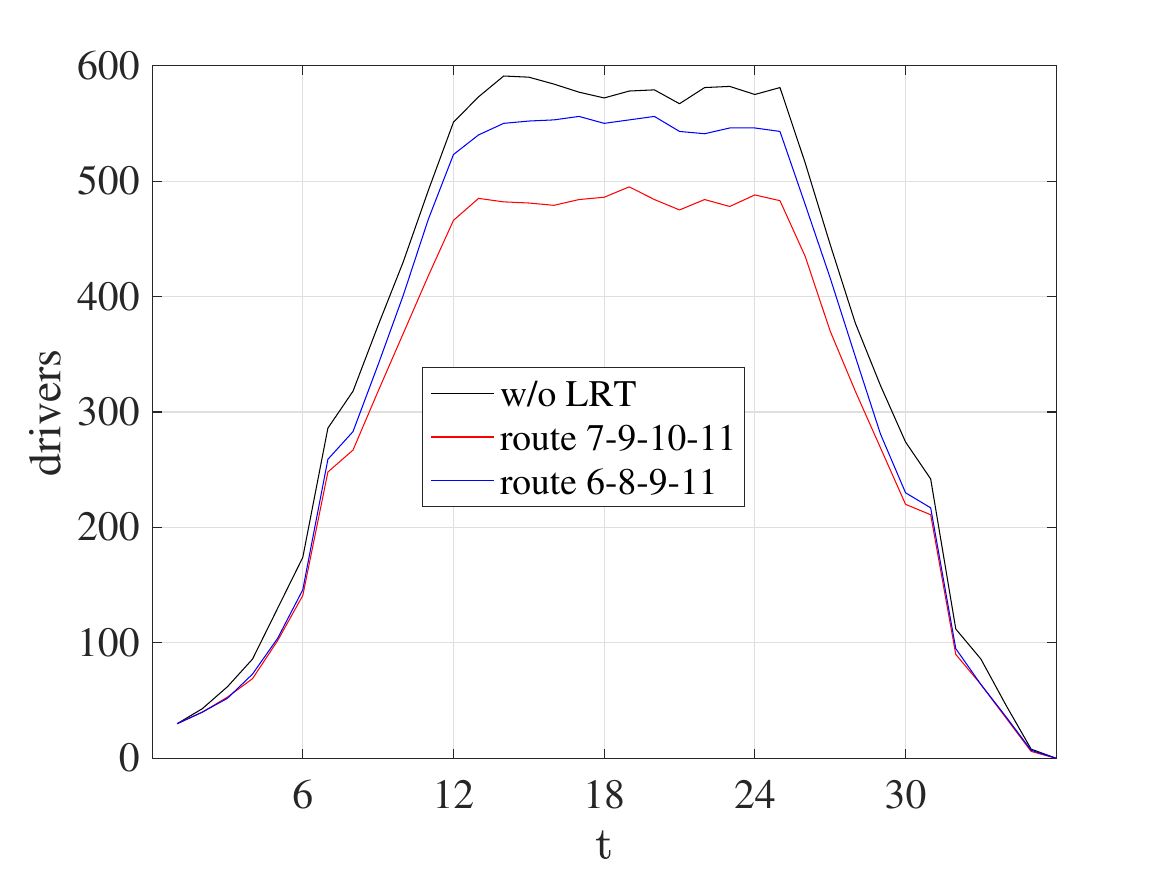}
    \caption{Number of drivers in downtown area \(\mathcal{D}=\{9, 10, 11\}\).}
    \label{fig: downtown}
\end{figure}

Notice that since drivers can start the game at different times during \(1\leq t\leq 12\) and they will only plan for the next 24 time steps. In other words, drivers with heterogeneous planning time windows will coexist in the network for \(t=2, 3, \ldots, 24\). Therefore the equilibrium of this game is a multi-commodity Markovian network equilibrium discussed in Section~\ref{subsection: models: multicommodity flow}. We consider the scenario where \(\alpha=10, \beta=0.2\), \(\gamma_{tss'}\) is given in Table~\ref{tab: gamma} and \(\text{dist}_{ss'}\) is given in Fig.~\ref{fig: seattle}. We compute the driver number in the downtown area \(\mathcal{D}=\{9, 10, 11\}\) by solving the optimization in \eqref{opt: multi optimal flow} using commercial software Mosek (\url{https://www.mosek.com}). The results are demonstrate in Fig.~\ref{fig: downtown}, where we can see that the driver number increases during \(1\leq t\leq 12\), then decreases during \(24\leq t\leq 36\). There are also two sudden changes in the increasing/decreasing rate around \(t=7\) and \(t=31\), which is due to the corresponding changes in values of \(\gamma\) in Table~\ref{tab: gamma}. This example extends the single commodity case considered in \cite{calderone2017markov} and \cite{li2019tolling} where all players enter and exit the game simultaneously. 

A relevant application of the above simulation framework is transportation network design. For example, suppose that Seattle city council is considering two candidate light rail transit (LRT) routes, 7-9-10-11 and 6-8-9-11 (see Fig.~\ref{fig: seattle}), as a means to alleviate the congestion caused by the ride-sharing traffic in downtown area, assuming that the LRT will reduce the demand of ride-sharing services (namely, value of \(\gamma_{tss'}\)) by \(50\%\) along its route. The simulated equilibrium with different LRT routes are also given in Fig.~\ref{fig: downtown}, which shows that route 6-8-9-11 is more effective that route 7-9-10-11 in terms of reducing amount of drivers in \(\mathcal{D}\). These results clearly demonstrate the power of Markovian network equilibrium model in transportation system design.

\begin{table}[ht]
\centering
\caption{Values of \(\gamma\) where \(\mathcal{D}=\{9, 10, 11\}\).}
\setlength{\tabcolsep}{3pt}
\begin{tabular}{ c|p{10mm}|p{10mm}|p{10mm}|p{10mm} } 
 \hline
 \(\gamma_{tss'}\) & \(s\notin\mathcal{D}\) \newline \(s'\in\mathcal{D}\)& \(s\in\mathcal{D}\) \newline \(s'\in\mathcal{D}\)& \(s\in\mathcal{D}\) \newline \(s'\notin\mathcal{D}\)& \(s\notin\mathcal{D}\) \newline \(s'\notin\mathcal{D}\) \\
 \hline
 \(1\leq t \leq 6\) & 600 & 200 & 60 & 60\\
 \(7\leq t \leq 30\)  & 200 & 400 & 200 & 60 \\
 \(31\leq t\leq 36\) & 60 & 100 & 600 & 60 \\
 \hline
\end{tabular}
\label{tab: gamma}
\end{table}

\subsection{Computation experiments}\label{subsection: Gurobi}

To demonstrate the efficiency of the algorithms developed in Section~\ref{section: algorithms}, we compare the computation time of our algorithms against commercial software Mosek, used in the previous section, over randomly generated examples. We use \(\mathrm{rand}(a, b)\) to denote a random number sampled from uniform distribution over interval \([a, b]\) where \(a, b\in\mathbb{R}\) and \(a\leq b\).
\begin{itemize}
    \item \(P_{sas'}=\mathrm{rand}(0, 1)\) for all \(s\in[S], a\in[A]\), then normalized  such that \(\sum_{s'}P_{sas'}=1\)
    \item \(\phi_{tsa}(\alpha)= \mathrm{rand}(1, 2)\alpha+\mathrm{rand}(1, 2)\) for all \(t\in[T], s\in[S], a\in[A]\).
    \item \(p_{ts}=\mathrm{rand}(0, 1)\) for all \(s\in[S]\) if \(t=1\) and zero otherwise.
\end{itemize}
In the variable demand case, we let \(\psi_{ts}(\alpha)=\mathrm{rand}(1, 2)\alpha-t+21\) for all \(t\in[T], s\in[S]\). In the multi-commodity flow case, we let \(\mathbb{T}=\{5, 10\}\).

\begin{figure}[htp]

\centering

\subfloat[Variable demand]{%
  \includegraphics[width=\linewidth]{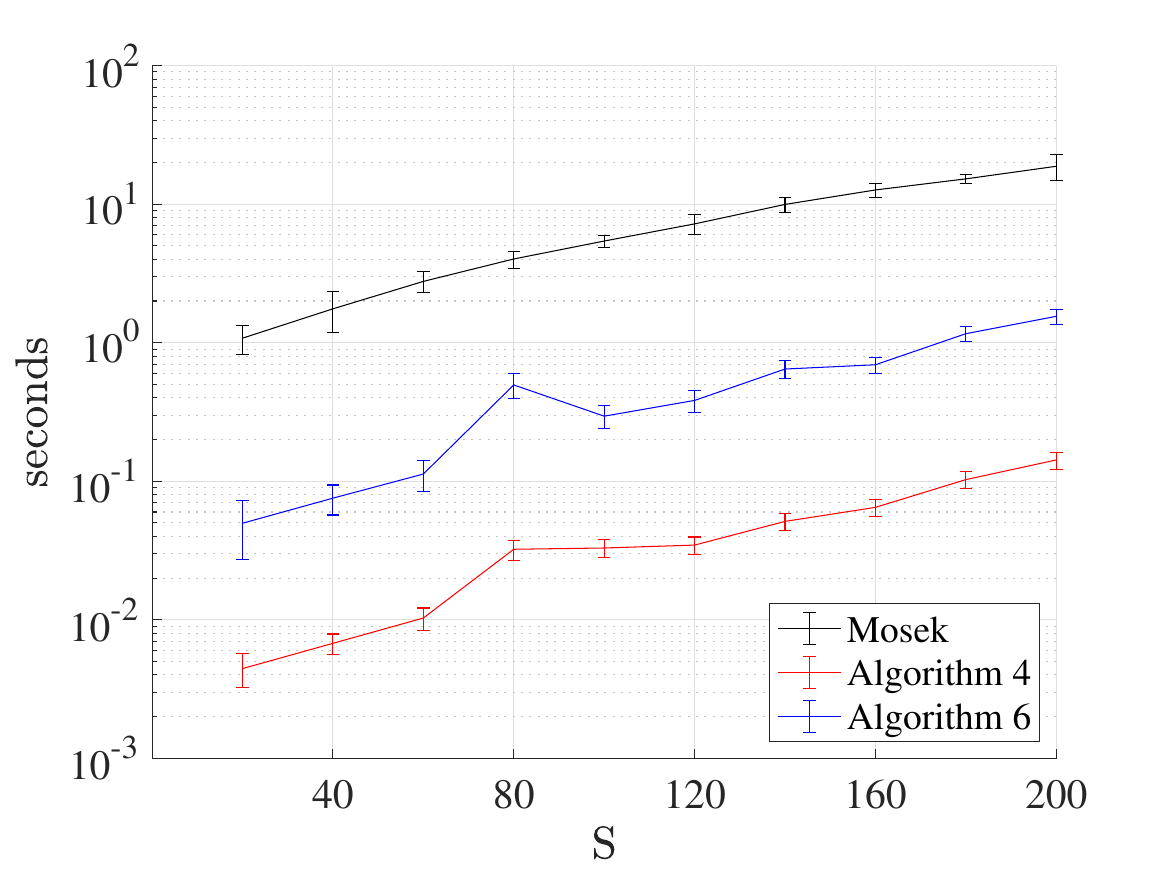} \label{fig: var}
}

\subfloat[Multi-commodity flow]{%
  \includegraphics[width=\linewidth]{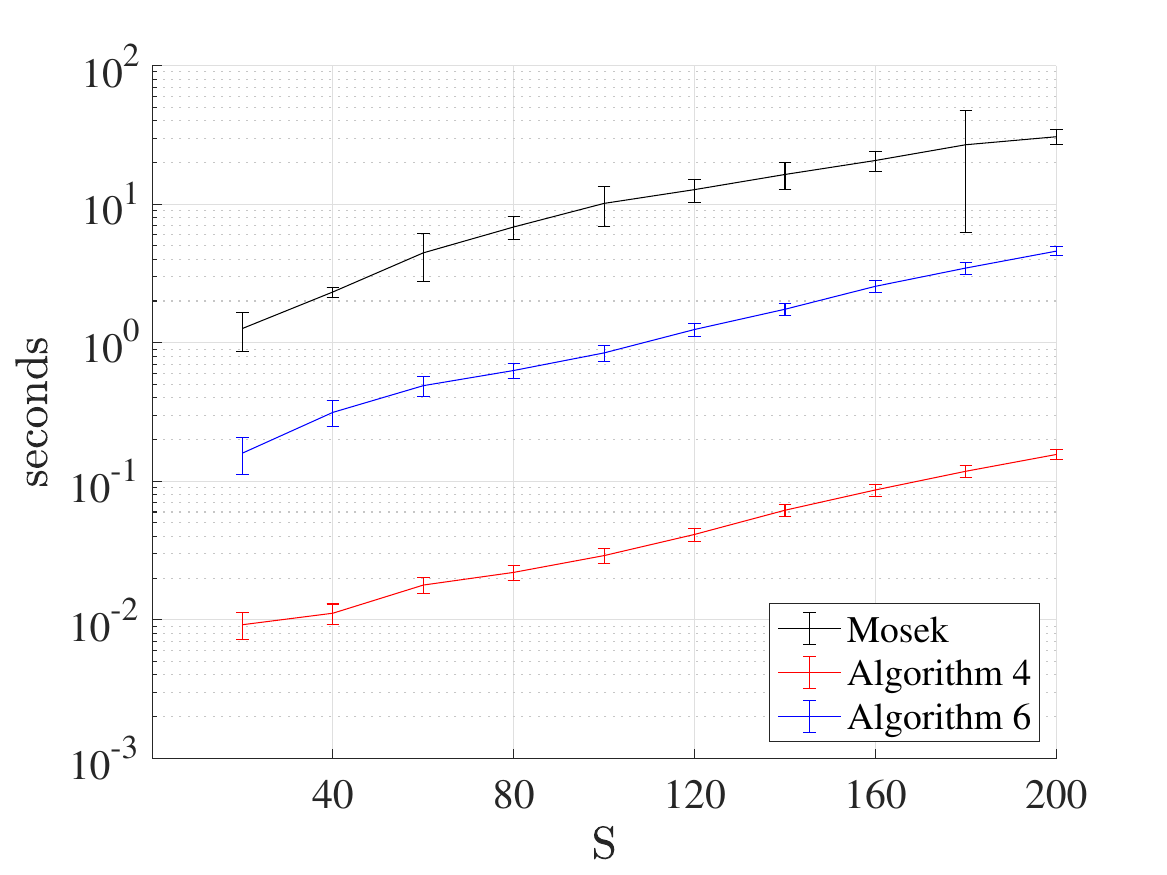} \label{fig: multi}
}

\caption{Average computation time and 3-standard deviation intervals of 100 experiments with \(T=A=10\).}
\label{fig: time}
\end{figure}

We fix \(T=A=10\) and let \(S\) range between \(20\) and \(200\), then test the computation time of Algorithm~\ref{alg: Frank-Wolfe}, Algorithm~\ref{alg: subgradient}, Algorithm~\ref{alg: multi Frank-Wolfe} and Algorithm~\ref{alg: multi subgradient}, where all algorithms terminate when their objective function value agrees with the optimal one obtained by Mosek with less than \(0.5\%\) relative error. The average computation time over 100 examples, along with corresponding 3-standard deviation interval are reported in Fig.~\ref{fig: time}. All codes are in MATLAB and run on a 1.6GHz laptop.
From results in Fig.~\ref{fig: time} we can see that, over the randomly generated 2000 examples, subgradient method and Frank-Wolfe method reduces the computation time consumed by Mosek by one and, respectively, two orders of magnitudes, at the price of a mere \(0.5\%\) of relative accuracy.

\section{Conclusion}\label{section: conclusion}
We study the variable demand and multi-commodity extensions in Markovian network equilibrium. We also propose efficient algorithms that outperform state-of-the-art commercial optimization software. However, the current work still has several limitations. For example, the cost of actions perceived by the players is assumed to be exact, rather than corrupted by stochastic noise, as considered in stochastic user equilibrium model. Further, the ending time of each player is fixed at the beginning of the game. A more realistic assumption is to allow the players to change their ending time and recompute the equilibrium periodically. We aim to address these limitations in future work.

\appendix 
\setcounter{equation}{0}
\renewcommand\theequation{A.\arabic{equation}}
\section{Appendix}
\subsection{Proof of Theorem~\ref{thm: variable demand}}\label{app: thm variable demand}
The objective function of problem \eqref{opt: optimal flow} is convex (since \(\phi_{tsa}\) and \(\psi_{ts}\) are strictly increasing), its constraints are affine, and the optimal value is obviously finite (since \(\phi_{tsa}\) and \(\psi_{ts}\) are finitely valued). These imply that a solution pair to \eqref{opt: optimal flow} and \eqref{opt: optimal potential} necessarily satisfy the KKT conditions \cite[Thm. 28.3.1]{rockafellar1970convex}. Let \(v_{ts}\) be the dual variable corresponding to the equality constraint containing \(p_{ts}\), let \(\mu_{tsa}, \theta_{ts}, \lambda_{ts}\geq 0\) be the dual variables corresponding to constraint \(y_{tsa}\geq 0\), \(z_{ts}\geq 0\) and, respectively, \(z_{ts}\leq p_{ts}\). Then the Lagrangian of \eqref{opt: optimal flow} and \eqref{opt: optimal potential} is given by
\begin{equation*}
\begin{aligned}
    &L(y, z,  v,  \mu, \lambda, \theta)= \textstyle\sum\limits_{t, s, a}\displaystyle\int_{0}^{y_{tsa}}\phi_{tsa}(\alpha)d\alpha-\textstyle\sum\limits_{t, s, a}\mu_{tsa}y_{tsa}\\
    &+\textstyle\sum\limits_{t,s}\displaystyle\int_{0}^{z_{ts}}\psi_{ts}(\alpha)d\alpha +\textstyle\sum\limits_{t,s,a} v_{ts}(p_{ts}-z_{ts}-y_{tsa})\\
    & \textstyle+\sum\limits_{t,t<T}\sum\limits_{s', a,s}v_{t+1,s}P_{s'as}y_{ts'a}-\textstyle\sum\limits_{t, s}((\theta_{ts}-\lambda_{ts})z_{ts}+\lambda_{ts}p_{ts}).
    \end{aligned}
\end{equation*}
The KKT conditions \cite[Thm.28.3]{rockafellar1970convex} of this Lagrangian include the following vanishing gradient conditions (by setting \(\partial L/\partial y_{tsa}, \partial L/\partial x_{ts}\) equal to zero)
\begin{equation}
        \begin{aligned}
        \textstyle v_{Ts}=&\phi_{Tsa}(y_{Tsa})-\mu_{Tsa},\\
        \textstyle v_{ts}=&\phi_{tsa}(y_{tsa})+\textstyle\sum\limits_{s'}P_{sas'}v_{t+1,s'}-\mu_{tsa},\enskip t\in[T-1],\\
        \textstyle v_{ts}=&\psi_{ts}(z_{ts})+\lambda_{ts}-\theta_{ts}, \enskip t\in[T],  \end{aligned}
        \label{KKT: vanishing gradient}
\end{equation}
for all \(s\in[S], a\in[A]\), and the complementarity conditions
\begin{equation}
    \begin{aligned}
         &y_{tsa}\mu_{tsa}=0, \enskip z_{ts}\theta_{ts}=0, \enskip \lambda_{ts}(z_{ts}-p_{ts})=0\\
         & y_{tsa}, \,z_{ts},\, \mu_{tsa}, \, \theta_{ts},\, \lambda_{ts}\geq 0, \enskip \forall t\in[T], s\in[S], a\in[A].
    \end{aligned}\label{KKT: complementarity}
\end{equation}
Combining \eqref{KKT: vanishing gradient} and \eqref{KKT: complementarity} yields \eqref{eqn: Wardrop quitting} and \eqref{eqn: Wardrop playing}. Note that same results can be derived from the dual problem \eqref{opt: optimal potential}.

\subsection{Proof of Theorem~\ref{thm: multicommodity} }\label{app: thm multicommodity}
The objective function of problem \eqref{opt: multi optimal flow} is convex (since \(\phi_{tsa}\) is increasing), its constraints are affine, and the optimal value is obviously finite (since \(\phi_{tsa}\) is finitely valued).These imply that a solution pair to \eqref{opt: multi optimal flow} and \eqref{opt: multi optimal potential} necessarily satisfies the KKT conditions \cite[Cor. 28.3.1]{rockafellar1970convex}. Let \(v^{\tau}_{ts}\) be the dual variable corresponding to the equality constraints containing \(p^{\tau}_ts\), let \(\mu^{\tau}_{tsa}\geq 0\) be the dual variables corresponding to constraint \(y^{\tau}_{tsa}\geq 0\). Then the Lagrangian of \eqref{opt: multi optimal flow} and \eqref{opt: multi optimal potential} is given by
\begin{equation*}
\begin{aligned}
     L(y, &\,v, \mu) =\textstyle\sum\limits_{t, s, a}\displaystyle\int_{0}^{\sum\limits_{\tau, \tau\geq t}y^{\tau}_{tsa}}\phi_{tsa}(\alpha)d\alpha\\
    &+\textstyle\sum\limits_{\tau}\sum\limits_{t,t\leq \tau}\sum\limits_{s, a} \left(v_{ts}^{\tau}(p_{ts}^{\tau}-y^{\tau}_{tsa})-\mu^{\tau}_{tsa}y^{\tau}_{tsa}\right)\\ &+\textstyle\sum\limits_{\tau}\sum\limits_{t,t<\tau}\sum\limits_{s', a,s}v^{\tau}_{t+1,s}P_{s'as}y^{\tau}_{ts'a}.
    \end{aligned}
\end{equation*}
The KKT conditions \cite[Thm.28.3]{rockafellar1970convex} of this Lagrangian include the following vanishing gradient conditions (by setting \(\partial L/\partial y^{\tau}_{tsa}\) equal to zero)
\begin{equation}
    \begin{aligned}
    \textstyle v^{\tau}_{\tau s}=&\phi_{\tau sa}\big(\textstyle\sum\limits_{j, j\geq \tau}y_{\tau sa}^{j}\big)-\mu^{\tau}_{\tau sa}=0,\\
    \textstyle v^{\tau}_{ts}=&\phi_{tsa}\big(\textstyle\sum\limits_{j, j\geq t}y^{j}_{tsa}\big)+\textstyle\sum\limits_{s'}P_{sas'}v^{\tau}_{t+1,s'}-\mu^{\tau}_{tsa}=0,
    \end{aligned}\label{KKT: multi vanishing gradient}
\end{equation}
for all \(t\in[\tau-1],\tau\in\mathbb{T}, s\in[S], a\in[A]\), and the complementarity conditions
\begin{equation}\label{KKT: multi complementarity}
    y_{tsa}^{\tau}\mu_{tsa}^{\tau}=0, \enskip y_{tsa}^{\tau},\, \mu_{tsa}^{\tau}\geq 0, 
\end{equation}
for all \(t\in[\tau],\tau\in\mathbb{T}, s\in[S], a\in[A]\). Combining \eqref{KKT: multi vanishing gradient} and \eqref{KKT: multi complementarity} yields \eqref{eqn: multi Wardrop}. Again, same results can be derived from the dual problem \eqref{opt: multi optimal potential}.
\subsection{Proof of Lemma~\ref{lem: variable demand}}\label{app: lem variable demand}
Using a similar argument as in the proof of Theorem~\ref{thm: variable demand}, we can show that the KKT conditions of \eqref{opt: linear optimal flow} are given by the following
\begin{equation}
    \begin{aligned}
    \textstyle\sum\limits_{a} y_{1sa}=&p_{1s}-z_{1s},\\
    \textstyle\sum\limits_{a}y_{t+1,sa}=&p_{t+1, s}-z_{t+1, s}+\textstyle\sum\limits_{s',a}P_{s'as}y_{ts'a},\, t\in [T-1],\\
     v_{Ts}=&u_{Tsa}-\mu_{Tsa},\\
        \textstyle v_{ts}=&u_{tsa}+\textstyle\sum\limits_{s'}P_{sas'}v_{t+1,s'}-\mu_{tsa},\enskip t\in[T-1],\\
        \textstyle v_{ts}=&w_{ts}+\lambda_{ts}-\theta_{ts},\\
         y_{tsa}\mu_{tsa}=&0, \enskip z_{ts}\theta_{ts}=0, \enskip \lambda_{ts}(z_{ts}-p_{ts})=0,\\
          y_{tsa},\, z_{ts}, \,&\mu_{tsa}, \, \theta_{ts},\, \lambda_{ts}\geq 0, 
    \end{aligned}\label{linear KKT: variable demand}
\end{equation}
for all \(t\in[T], s\in[S], a\in[A]\). Let \((\hat{v}, \hat{\pi})\) be the output of Algorithm~\ref{alg: Bellman} with input \((P, u, T)\), \(\hat{z}=p\odot (\hat{v}>w)\), \(\hat{y}\) be the output of Algorithm~\ref{alg: Kolmogorov} with input \((\hat{\pi}, p-\hat{z}, P, T)\), let 
\begin{equation}
    \begin{aligned}
    \hat{\mu}_{Tsa}=&-\hat{v}_{Ts}+u_{Tsa},\\
        \textstyle \hat{\mu}_{tsa}=&-\hat{v}_{ts}+u_{tsa}+\textstyle\sum\limits_{s'}P_{sas'}\hat{v}_{t+1,s'},\enskip t\in[T-1],\\
        \hat{\lambda}_{ts}=&\max\{\hat{v}_{ts}-w_{ts}, 0\}, \enskip \hat{\theta}_{ts}=\max\{w_{ts}-\hat{v}_{ts}, 0\},
    \end{aligned}
\end{equation}
for all \(t\in[T], s\in[S], a\in[A]\). Then it is straightforward to verify that \((\hat{y}, \hat{z}, \hat{v},\hat{\mu}, \hat{\lambda},\hat{\theta})\) satisfies all the KKT conditions in \eqref{linear KKT: variable demand}, hence \((\hat{y}, \hat{z})\) solves \eqref{opt: linear optimal flow}, which proves the equality. The inequality follows from the fact that, when \((u, w)\) in \eqref{opt: linear optimal flow} is perturbed to another value \((u', w')\), solution \((\hat{y}, \hat{z})\) is still feasible, but can be suboptimal.

\subsection{Proof of Lemma~\ref{lem: multicommodity}}\label{app: lem multicommodity} Notice that in optimization \eqref{opt: linear multi optimal flow}, both objective function and constraints are completely separable across \(y^\tau\) with different value of \(\tau\). In other words, solving \eqref{opt: linear multi optimal flow} is equivalent to solve the following optimization problem for each value of \(\tau\in\mathbb{T}\) separately
\begin{equation}
    \begin{array}{ll}
        \underset{y^\tau}{\mbox{min.}} & \sum\limits_{t\leq \tau,s,a} u_{tsa}y_{tsa}^\tau \\
         \mbox{s.t.} & \sum\limits_{a} y^{\tau}_{1sa}=p^{\tau}_{1s}\\ &\sum\limits_{a}y^{\tau}_{t+1, sa}=p^{\tau}_{t+1,s}+\sum\limits_{s',a}P_{s'as}y^{\tau}_{ts'a},\enskip t\in[\tau-1],\\
         &0\leq y^{\tau}_{tsa}, \enskip
         \forall t\in[\tau], s\in[S], a\in[A].
    \end{array}
    \label{opt: linear single optimal flow}
\end{equation}
Since problem \eqref{opt: linear single optimal flow} is nothing but an instance of \eqref{opt: linear optimal distribution} with \(T=\tau\), it can be solved by the output of Algorithm~\ref{alg: Kolmogorov} with input \((\hat{\pi}^\tau, p^\tau, P, \tau)\), where \(\hat{\pi}^\tau\) is the output of Algorithm~\ref{alg: Bellman} with input \((P, u, \tau)\). This proves the equality; the inequality follows from the fact that, when \(u\) in \eqref{opt: linear multi optimal flow} is perturbed to another value \(u'\), solution \((\hat{y}^\tau, \tau\in\mathbb{T})\) is still feasible, but can be suboptimal.

\subsection{Proof of Theorem~\ref{thm: FW convergence} }\label{app: thm FW convergence}

We start with Algorithm~\ref{alg: Frank-Wolfe}. The per-iteration computation of Algorithm~\ref{alg: Frank-Wolfe} is clearly dominated by the execution of Algorithm~\ref{alg: Bellman} and Algorithm~\ref{alg: Kolmogorov}, which together cost \(O(\sigma TS^2A)\) arithmetical operations. Let \(f(y, z)\) denote the objective function of problem~\eqref{opt: optimal flow}. Then the gradients of \(f\) is given by
\[\partial_y f=\phi(y), \enskip \partial_z f=\psi(z), \]
where \(\phi(y)\) and \(\psi(z)\) are defined as in \eqref{eqn: compact tensor}. Then under Assumption~\ref{asp: Lipschitz variable demand}, we know both \(\frac{\partial f}{\partial y}\) and \(\frac{\partial f}{\partial z}\) are Lipschitz. In addition, for any \((y, z)\) satisfying the constraints of problem \eqref{opt: optimal flow}, one must have \(y_{tsa}\in[0, \rho]\) and \(z_{ts}\in[0, \rho]\), due to Assumption~\ref{asp: variable demand}. In other words, the constraint set of problem \eqref{opt: optimal flow} is a subset of \([0, \rho]^{T\times S\times A}\times [0, \rho]^{T\times S}\), which is bounded.

Therefore problem \eqref{opt: optimal flow} is minimizing a function with Lipschitz gradients over a bounded set. Hence the Frank-Wolfe method given by Algorithm~\ref{alg: Frank-Wolfe} converges to \(\epsilon\)-optimal solution in \(O(\frac{1}{\epsilon})\) iterations \cite[Thm.3.8]{bubeck2015convex}. The proof for Algorithm~\ref{alg: multi Frank-Wolfe} is similar.

\subsection{Proof of Theorem~\ref{thm: subgrad convergence} }\label{app: subgrad convergence}
We start with Algorithm~\ref{alg: subgradient}. First, the per-iteration computation of Algorithm~\ref{alg: subgradient} is clearly dominated by the execution of Algorithm~\ref{alg: Bellman} and Algorithm~\ref{alg: Bellman}, which together cost \(O(\sigma TS^2A)\) arithmetical operations. From Assumption~\ref{asp: variable demand} and \ref{asp: Lipschitz variable demand} we have, for any \(u_{tsa}, u'_{tsa}, w_{ts}, w'_{ts}\in[0, \rho]\)
\begin{equation*}
    \begin{aligned}
    |u_{tsa}-u'_{tsa}|=&|\phi_{tsa}(\phi_{tsa}^{-1})(u_{tsa})-\phi_{tsa}(\phi^{-1}_{tsa})(u'_{tsa})|\\
    \leq & L|\phi_{tsa}^{-1}(u_{tsa})-\phi_{tsa}(\phi^{-1}_{tsa})(u'_{tsa})|,\\
    |w_{ts}-w'_{ts}|=&|\psi_{ts}(\psi_{ts}^{-1})(w_{ts})-\psi_{ts}(\psi^{-1}_{ts})(w'_{ts})|\\
    \leq & L|\phi_{tsa}^{-1}(u_{tsa})-\phi_{tsa}(\phi^{-1}_{tsa})(u'_{tsa})|,
    \end{aligned}
\end{equation*}
for all \(t\in[T], s\in[S], a\in[A]\), which implies the objective function of problem \eqref{opt: outer optimal potential} (in particular, the integral terms) is \(\frac{1}{L}\)-strongly convex \cite[Thm.2.1.10]{nesterov2013introductory}. Let \(-f(u, w)\) denote the objective function of problem \eqref{opt: outer optimal potential}. Then from Lemma~\ref{lem: variable demand} we know that the subgradients of \(f\) is given by 
\[\partial_u f=-\hat{y}+\phi^{-1}(u), \enskip \partial_w f=-\hat{z}+\psi^{-1}(w),\]
where \(\phi^{-1}(u), \psi^{-1}(w)\) are defined as in \eqref{eqn: compact tensor}, \((\hat{y}, \hat{z})\) is a solution to problem \eqref{opt: linear optimal flow}. From Assumption~\ref{asp: variable demand} we know that \(\phi^{-1}(u)\in[0, \rho]^{T\times S\times A}\) and \(\psi^{-1}(w)\in[0, \rho]^{T\times S}\). Further, \((\hat{y}, \hat{z})\) must satisfy the constraints in problem \eqref{opt: linear optimal flow}, which implies that \(\hat{y}\in[0, \rho]^{T\times S\times A}\) and \(\hat{z}\in[0, \rho]^{T\times S}\). Hence the elements in \(\partial_u f\) and \(\partial_w f\) are bounded. 

Therefore, problem \eqref{opt: outer optimal potential} is minimizing a strongly convex function whose subgradients have bounded elements. hence the projected subgradient method given by Algorithm~\ref{alg: subgradient} converges to an \(\epsilon\)-optimal solution in \(O(\frac{1}{\epsilon})\) iterations \cite[Th,.3.9]{bubeck2015convex}. The proof for Algorithm~\ref{alg: multi subgradient} is similar.

\bibliographystyle{apalike}      
\bibliography{reference}           

\end{document}